\documentclass[10pt]{article}
\usepackage{amsmath,amsfonts,amsthm,amssymb,amscd,authblk}
\numberwithin{equation}{section}

\newcommand{\be}{\begin{equation}}
\newcommand{\ee}{\end{equation}}
\newcommand{\bs}{\begin{split}}
\newcommand{\es}{\end{split}}
\newcommand{\ba}{\begin{align}}
\newcommand{\ea}{\end{align}}
\newcommand{\basl}[1]{\begin{align}\begin{split}\label{#1}}
\newcommand{\bas}{\begin{align}\begin{split}}

\newtheorem{theo}{Theorem}[section]
\newtheorem{prop}[theo]{Proposition}
\newtheorem{lemm}[theo]{Lemma}

\newcommand\fpr{\hfill$\Box$\null}

\newcommand\N{\mathbb{N}}

\newcommand\R{\mathbb{R}}
\newcommand\C{\mathbb{C}}

\title{On the anti-Wick symbol as a Gelfand-Shilov generalized function}

\begin{document}

\author[1]{L. Amour}
\author[2]{N. Lerner}
\author[1]{J. Nourrigat}

\affil[1]{{\sc LMR FRE CNRS 2011}, Universit\'e de Reims, France.}
 \affil[2]{ {\sc IMJ UMR CNRS 7586}, Sorbonne Universit\'e, France.}

\date{}

\maketitle

\vskip 0.3cm

\begin{abstract}\noindent
The purpose of this article is to prove that the anti-Wick symbol of an operator mapping $ {\cal S}(\R^n)$ into ${\cal S}'(\R^n)$, which is generally not a tempered distribution,  can still be defined as a Gelfand-Shilov generalized function. This result relies on test function spaces embeddings involving the Schwartz and Gelfand-Shilov spaces. An additional embedding concerning Schwartz and Gevrey spaces is also given.
\end{abstract}

\

\noindent
{\it Keywords:} Anti-Wick symbol, Gelfand-Shilov generalized functions, pseudodifferential calculus, Gevrey spaces, test function spaces embeddings.

\

\noindent
{\it MSC 2010:} 47G30,46F05

\parindent=0pt
\parindent = 0 cm

\parskip 10pt
\baselineskip 12pt

 \section{Statement of the result.}    \label{s1}

For any operator $A$ mapping ${\cal S}(\R^n)$ into ${\cal S}'(\R^n) $, it is well known that the
  Weyl symbol $\sigma ^{Weyl } (A)  $ and the distribution
kernel $K_A$ are tempered distributions on $\R^{2n}$   satisfying
(with the notations of \cite{Ler}):
\be\label{def-Weyl-1}
K_A(x , y)  =  \int _{\R^n} e^{2 i \pi  ( x-y) \cdot \xi }
 \sigma ^{Weyl } (A)  \left ( \frac {x +  y} {2}  , \xi \right ) d\xi,\ee
or equivalently:
\be\label{def-Weyl-2}   \sigma ^{Weyl } (A) (x , \xi ) =
 \int _{\R^n} e^{ - 2i\pi  t \cdot \xi }
K_A  \left( x + \frac{t}{2}  , x - \frac{t}{2}  \right ) dt.\ee

(The above spaces  ${\cal S}(\R^n)$ and ${\cal S}'(\R^n) $ refer respectively to the standard  Schwartz space and its topological dual).

The  two identities (\ref{def-Weyl-1}) and (\ref{def-Weyl-2}) are naturally  understood in the sense of tempered distributions. Recall that these equalities can be regarded as a basis of Weyl pseudodifferential  calculus \cite{HO}.

Besides, an operator $A$ can also be given under the form:
 \be\label{anti-Wick}  <Af , g> = (2\pi )^{-n} \int _{\R^{2n}} F(X)
 < f, \Psi_X> < \Psi _X , g> dX \ee
 where the $\Psi _X $ are the coherent states:
 \be\label{CS} \Psi _X  = \tau _X \Psi_0, \hskip 0.5cm
 \Psi_0 (u)= 2 ^\frac{n}{4}e^{ - \pi |u|^2}, \hskip 0.5cm
 \tau _ {(x , \xi)} f (u) = f(u-x) e^{2 i\pi \left( u - \frac{x}{2} \right) \cdot \xi }
 \hskip 0.5cm  (X=(x,\xi)\in\R^{2n},\ u\in\R^n) \ee
and where $F$ is, for instance,  a bounded continuous function on $\R^{2n}$, called
 anti-Wick symbol of $A$. This symbol can be denoted  $\sigma ^{AW} (A)$.
For that purpose, see, {\it e.g.},  \cite{FAB,C-R,F,Ler}.

In this case, it is well known that:
 \be\label{heat-2} \sigma ^{Weyl} (A) = 2^n  \sigma ^{AW} (A) \ast e^{ - 2 \pi \Gamma  }
 \quad\text{where}\quad \Gamma (X)= |X|^2  \ee
 (where $\ast$ is the standard convolution product)
 or equivalently:
 \be\label{heat-2+} \sigma ^{Weyl} (A) = e^{\frac {1} {8 \pi}  \Delta }  \sigma ^{AW} (A) \ee
(where $\Delta$ denotes  the Laplacian operator).

Our aim here is to
 prove that it is always possible to define the anti-Wick
 symbol of any operator $A$ mapping ${\cal S}(\R^n)$ into ${\cal S}'(\R^n) $.
This symbol is not a tempered distribution, but a Gelfand-Shilov  generalized function
(\cite{GS-1,GS-2}).

First, one can a priori define a linear form   $T(A) $ on the space
$e^{ \frac {1} {8 \pi}  \Delta } {\cal S} (\R^n)$ by:
$$ < T(A) , e^{ \frac {1} {8 \pi}   \Delta } \varphi > = < \sigma ^{Weyl} (A) ,
\varphi >, \quad \varphi \in {\cal S} (\R^n).$$
One then has:
   \be\label{heat-3}  \sigma ^{Weyl} (A) = e^{ \frac {1} {8 \pi}   \Delta }   T_A.  \ee
As a consequence and by analogy with (\ref{heat-2}), one can consider  that $T(A)$ is by definition the anti-Wick symbol of
 $A$. 
 
 Next, Theorem \ref{t1} below proves that  $T(A)$ defined above
is actually also a Gelfand-Shilov generalized function.

Before stating Theorem \ref{t1},  let us recall here the definition of the space $S(\lambda , \mu) (\R^n)$ of these test functions used by Gelfand and Shilov \cite{GS-1,GS-2}. See also, $e.g.$ \cite{Ler-GS1,Ler-GS2,Ler-GS3} for applications and \cite{Smi} for related spaces.

The space $S(\lambda , \mu) (\R^n)$ refers to the space of functions
  $\varphi \in C^{\infty}(\R^n)$ such that, there exists a constant $A>0$ satisfying
for all multi-indices $\alpha $ and $\beta $, for all $x\in \R^n$:
   \be\label{fct-test} |x^{\alpha } \partial^{\beta }_x \varphi (x) | \leq A^{|\alpha |+ |\beta  | }
   ( \alpha !) ^{\lambda } ( \beta !) ^{\mu } . \ee

The following Theorem is proved in Section \ref{s2}.

   \begin{theo} \label{t1} If $\lambda >0 $ and  $0 < \mu < \frac{1}{2}$ then
the space $S(\lambda , \mu) (\R^n)$ is continuously embedded in the space
   $e^{ \frac {1} {8 \pi}   \Delta } {\cal S} (\R^n)$.

    \end{theo}

  As a consequence, one obtains that the anti-Wick symbol $T(A)$ of any operator 
   $A$ mapping $ {\cal S}(\R^n)$ into ${\cal S}'(\R^n) $,   is well defined (by restriction) as a continuous linear form on  $S(\lambda , \mu) (\R^{2n})$ for any
 $\lambda >0 $ and any $\mu < \frac{1}{2}$. That is,
 $T(A)$ is a Gelfand-Shilov generalized function.

 Let us also mention the following fact as a  complementary result  concerning anti-Wick symbols. In  \cite{A-J}, one provides conditions written in terms of the action  of an operator 
$A$ on coherent states to ensure that the anti-Wick symbol  of $A$ is a bounded continuous function on $\R^{2n}$.
Also note that this latter result is actually to be compared with Unterberger result \cite{U} giving a similar necessary and sufficient condition in order that the Weyl symbol of $A$ is a 
 $C^{\infty}$ function on $\R^{2n}$, being bounded together with all of its derivatives.
 
 Theorem \ref{t1} is proved in Section \ref{s2} and a related test function spaces embedding concerning Gevrey spaces is given in Section \ref{s3}. 
 
   \section{Proof of Theorem \ref{t1}.}    \label{s2}

The Proposition below follows from Proposition  6.1.8 of Nicola-Rodino \cite{NR} but we give a proof for the sake of the reader convenience.

    \begin{prop} \label{pGS-Holo}
    Suppose that $\lambda >0$ and $0 < \mu < 1$. Then, any function in  $S(\lambda , \mu) (\R^n)$
  extends to a holomorphic function $u$ on  $\C^n$ satisfying for some constant $K>0$:
  \be\label{majo-prol-holo}    e^{ \varphi (x) }    |u(x+i y) | \leq K e^{\psi (y) }, \ee
  with
   \be\label{phi-psi}  \varphi (x) =  \frac {\lambda  } {2} \sum _{j=1}^n \left|\frac{ x_j}{ A}\right |^{\frac{1}{\lambda} }\quad\text{and}\quad
  \psi (y) = 2 (1 - \mu)
   \sum _{j=1}^n | A  y_j  |^{\frac{1}{1 - \mu} },   \ee
where the constant $A>0$ is the constant in (\ref{fct-test}).

    \end{prop}

    {\it Proof of Proposition \ref{pGS-Holo}.} 
     For any $x>0$, one has:
  $$ \sup _{k \geq 0}  \frac {x^k} {k!}
  \geq \frac {1} {2} e^{ \frac{x}{2}}.$$
If $\lambda >0$, this yields for all $x\in \R^n$:
  $$ \frac {1} {2^{n \lambda } } e^{ \varphi (x) } \leq  \sup _{\alpha }
  \frac { |x^{\alpha } | } { A^{\alpha } (\alpha ! ) ^{\lambda } }.
    $$
If $\nu >0 $ and $x>0$, one has:
   $$  \left ( \frac {x^k} {k!} \right )^{\nu } \leq
     \frac {1} {2^{k \nu} } e^{2 \nu x} $$
implying
  $$ \sum _{k=0}^{\infty } \left ( \frac {x^k} {k!} \right )^{\nu }
  \leq  C e^{2 \nu   x}\quad\text{with}\quad C = \sum _{k=0}^{\infty }  \frac {1} {2^{k \nu} }. $$
Therefore, for all $u\in S(\lambda , \mu) (\R^n)$ with $\mu < 1$ and
 $\nu = 1 - \mu $:
  $$ \sum _{\beta } \left | \partial ^{\beta }_x u(x) \frac { (i y )^{\beta } }
  {\beta ! } \right | \leq  C^n e^{\psi (y) }.
   $$
Thus inequality (\ref{majo-prol-holo}) with (\ref{phi-psi}) is valid and the proof is then complete.

  \fpr

Next, we denote by  $E(\C^n)$ the space of holomorphic functions $\varphi$
on $\C^n$ satisfying for all $m\geq 0$:
 \be\label{space}
 \int _{{\bf C }  ^{n } } e^{ - 2 \pi  |{\rm Im} z|^2 }
 ( 1 + |{\rm Re } z| )^m   |\varphi (z) | dz < \infty.   \ee

  \begin{prop}\label{stepA}  If $\lambda >0 $ and $\mu < \frac{1}{2}$ then any element of 
  $S(\lambda , \mu) (\R^n)$
   can be extended to a holomorphic function belonging to $E(\C^n)$. This defines a continuous embedding of $S(\lambda , \mu) (\R^n)$ into $E(\C^n)$.

     \end{prop}

   {\it Proof of Proposition \ref{stepA}.} 
   One notices for the function $\psi$ defined in (\ref{phi-psi}) that: 
   $$ e^{\psi (y)} =  e^{ 2 \pi  |y|^2 } g(y) $$
where the function $g$ belongs to $L^1 (\R^n)$ if $\mu < \frac{1}{2}$. Besides, the function $e^{-\varphi (x)}$ is rapidly decreasing if $\lambda >0$. Then, Proposition \ref{stepA} follows from Proposition \ref{pGS-Holo}.

\fpr

The following Fourier transform is used in the sequel:
$$ \widehat { u } (\xi )  =  \int _{\R^n}
 u (x) e^{-2 i \pi x\cdot \xi } dx.
 $$

 \begin{prop}\label{stepB}   For any  $u\in E(\C^n)$, there is a function
 $\Phi\in {\cal S}(\R^n)$ such that $ e^{ \frac {1} {8 \pi } \Delta } \Phi $ 
is equal to $u$ restricted to $\R^n$.

 \end{prop}

  {\it Proof of Proposition \ref{stepB}.} 
   Let $u$ be any function belonging to $E(\C^n)$. Then,  one has:
  $$ \widehat { u } (\xi )  =  \int _{\R^n}
  u (x) e^{-2 i \pi x\cdot \xi } dx  =  \int _{\R^n}
  u  (x+ i y ) e^{-2 i\pi ( x+ i y) \cdot \xi } dx,$$
  for any $ \xi \in \R^n$. 
This implies:
  $$ e^{-  2 \pi ( |y|^2 + y \cdot \xi )  }  \widehat { u } (\xi )   =
  \int _{\R^n}
 u (x+ i y ) e^{- 2 \pi ( |y|^2    + i x \cdot \xi ) }dx. $$
 Therefore, one deduces integrating with respect to $y$:
 $$  2 ^{- \frac{n}{2}}  e^{ \frac {\pi } {2} |\xi |^2 }    \widehat {  u } (\xi )  =   \int _{\C^n}
 u (x + i y) e^{- 2 \pi ( |y|^2  + i x \cdot \xi )  } dx dy. $$
Set:
  $$ \kappa (\xi) = 2 ^{ \frac{n}{2}}   \int _{\C^n} u (x + i y)
  e^{- 2 \pi ( |y|^2  + i x \cdot \xi )  } dx dy. $$
This function $\psi$ belongs to ${\cal S} (\R^n)$ and verifies:
 $$  e^{- \frac {\pi } {2}  |\xi |^2} \kappa (\xi) =  \widehat u(\xi).   $$
 Thus, one observes that the function  $\Phi$ defined by:
$$ \Phi  (x )  =  \int _{\R^n}
 \kappa (\xi) e^{2 i \pi x\cdot \xi } d\xi
 $$
is indeed satisfying:
  $$ e^{ \frac {1 } {8 \pi } \Delta } \Phi = u.$$
  
  \fpr

{\it Proof of Theorem \ref{t1}.}  It is a direct consequence of Propositions \ref{stepA} et \ref{stepB} when replacing  $n$ by $2n$.

\fpr

 \section{An additional embedding.}    \label{s3} 
 
 We give here a supplementary embedding involving Schwartz and Gevrey spaces.
 
 We begin with the following Lemma.
 
 \begin{lemm} \label{der-Gauss}One has: 
$$ \left |  \frac  {d^m} {d x^m} e^{ - \frac {x^2}   {2} } \right  |
 \leq \sqrt {2} (2 \pi ) \frac{1}{4} \sqrt {m! }   (m+1 ) ^\frac{1}{4}, $$
for any $m \geq 0$ and all $x\in \R$:

\end{lemm}

{\it Proof of Lemma \ref{der-Gauss}.} 
Set:
$$ f_m (x) = \frac  {d^m} {d x^m} e^{ - \frac {x^2}   {2} } = (-1)^m
e^{ - \frac {x^2}   {2} } H_m (x) $$
where the $H_m$ are the Hermite functions.

One checks that:
$$ | f_m (x)|^2 \leq 2 \int _ {\R} | f_m (t)| | f'_m (t)| dt
\leq 2 \Vert f_m \Vert \  \Vert f_{m+1}  \Vert $$
where the above norms are the $L^2(\R)$ norms. Besides, one has:
$$  \Vert f_m \Vert ^2 = \int _{\R}  e^{ - x^2 }  |H_m (x) |^2 dx
\leq \int _{\R}  e^{ - \frac {x^2}   {2} }  |H_m (x) |^2 dx =
\sqrt {2 \pi } m!. $$
The proof of the Lemma thus follows.

\fpr

The Gevrey type space $G^s$ considered here is defined  as the space of smooth functions $f$ such that, there exist $K>0$ and $C>0$ satisfying for all multi-indices $\alpha\in \N^n$ and all $x\in\R^n$, 
$|\partial_x^\alpha f(x)|\leq K C^{|\alpha|}(|\alpha| !)^s$.

We then obtain the  next result.

 \begin{prop} \label{inclusion-supp}  If the function $f\in {\cal S}(\R^n)$ then the function
$ e^{ \frac {1} {8 \pi } \Delta } f $ 
belongs to the Gevrey space  $G^{\frac{1}{2}} (\R^n)$.

 \end{prop}

 {\it Proof of Proposition  \ref{inclusion-supp}.}  
Recall that:
 $$ e^{ \frac {1} {8 \pi } \Delta } f = 2^{ \frac{n}{2}}   f \ast e^{ - 2 \pi \Gamma  },
 \quad \Gamma (x)= |x|^2   $$
for any   $f\in {\cal S } (\R^n)$.
According to Lemma  \ref{der-Gauss} above, there exists $C>0$ satisfying:
 $$ \left |  \partial_x ^{\alpha }  e^{ - 2 \pi \Gamma (x) }     \right |
  \leq C ^{1 + |\alpha | } (\alpha !) ^{\frac{1}{2}} $$
for any multi-indice $\alpha\in \N^n$ and all $x\in\R^n$.

This proves Proposition  \ref{inclusion-supp}.

\fpr

laurent.amour@univ-reims.fr\newline
LMR, Universit\'e de Reims Champagne-Ardenne,
 Moulin de la Housse, BP 1039,
 51687 REIMS Cedex 2, France.
 
nicolas.lerner@imj-prg.fr\newline
IMJ-PRG, Sorbonne Universit\'e, 
Campus Pierre et Marie Curie,
4 place Jussieu, 75252 Paris Cedex, France.

jean.nourrigat@univ-reims.fr\newline
LMR, Universit\'e de Reims Champagne-Ardenne,
 Moulin de la Housse, BP 1039,
 51687 REIMS Cedex 2, France.

 \end{document}